\documentclass[11pt]{article}
\usepackage{amssymb}
 \csname @addtoreset\endcsname{equation}{section}
\newtheorem{theorem}{Theorem}
\newtheorem{lemma}{Lemma}
\newtheorem{definition}{Definition}

 \def\e{\varepsilon}

\def\defi{\stackrel{{\scriptscriptstyle \Delta}}{=}}

\def\a{\alpha}
\def\d{\delta}

\def\w{\widehat}

\def\esssup{\mathop{\rm ess\, sup}}
\def\essinf{\mathop{\rm ess\, inf}}

\def\R{{\bf R}}

\def\L{L}

\def\g{\gamma}

\def\W{{\cal W}^*}

\def\oo{\bar}

\def\p{\partial}

\def\V{{\cal V}}

\def\L{{\cal L}}

\newcommand{\be}{\begin{equation}}
\newcommand{\ee}{\end{equation}}
\newcommand{\bd}{\begin{displaymath}}
\newcommand{\ed}{\end{displaymath}}
\newcommand{\ba}{\begin{array}{ll}}
\newcommand{\ea}{\end{array}}
\newcommand{\baa}{\begin{eqnarray}}
\newcommand{\eaa}{\end{eqnarray}}
\newcommand{\baaa}{\begin{eqnarray*}}
\newcommand{\eaaa}{\end{eqnarray*}}
\font\sm=cmr10
\def\L{{\cal L}}

\def\W{{\cal W}}

\def\Q{{\cal Q}}

\def\R{{\bf R}}

\date{Submitted: September 27, 2010. Revised: February 17, 2011  }
\title
{ On prescribed  change of  profile for solutions of parabolic
equations  \footnote{{\it Journal of Physics A: Mathematical and
Theoretical} (May 2011) {\bf 44} 225204.}}
\author{
Nikolai Dokuchaev\\  {\sm Department of Mathematics \& Statistics,
Curtin University,}\\
{\sm  GPO Box U1987, Perth, 6845 Western Australia}\\ {\sm email
N.Dokuchaev@curtin.edu.au}}
\begin{document}
\maketitle
\begin{abstract}
Parabolic equations with homogeneous Dirichlet conditions  on the
boundary are studied in a setting where the solutions are required
to have a prescribed change of the profile in fixed time, instead of
a Cauchy condition. It is shown that this problem is well-posed in
$L_2$-setting. Existence and regularity results are established, as
well as an analog of the maximum principle.
\\
PACS 2010: 02.30.Zz,    
02.30.Jr.   
\\
{\it Abbreviated title:} On prescribed  absorption for parabolic
diffusions
\end{abstract}
\section{Introduction}
Parabolic diffusion equations have fundamental significance for
natural and social sciences, and various boundary value problems for
them were widely studied including   inverse and ill-posed problems;
see examples in Miller (1973),  Tikhonov and Arsenin (1977), Glasko
(1984), Prilepko {\it et al} (1984),
  Beck (1985), Seidman (1996).
According to Hadamard criterion, a boundary value problem is
well-posed if there is existence and uniqueness of the solution, and
if there is continuous dependence of the solution on the boundary
data. Otherwise, a problem is ill-posed.
\par
Apparently there are boundary value problems that do not fit the
framework given by the classical theory of well-posedness (see
examples in Dokuchaev (2007,2010)).
\par
For  parabolic equations, it is commonly recognized that the type of
the boundary conditions usually defines if a problem is well-posed
or ill-posed.
  A classical example is the
heat equation \baaa \frac{\p u}{\p t}=\frac{\p^2 u}{\p x^2}, \quad
t\in[0,T]. \label{sample}\eaaa The problem for this equation with
the Cauchy condition at initial time $t=0$ is well-posed in usual
classes of solutions, including classical, H\"older and square
integrable solutions.  In contrast, the problem with the Cauchy
condition at terminal time $t=T$  is ill-posed for this heat
equation for all these classes. In particular, this means  that a
prescribed profile of temperature at time $t=T$ cannot be achieved
via an appropriate selection of the initial temperature. In
addition, $L_2$-norms of solutions cannot be estimated by
$L_2$-norms of the boundary data (i.e, the dependence on boundary
data  is not continuous). This makes this problem ill-posed, despite
the fact that solvability and uniqueness still can be achieved for
some very smooth analytical boundary data or for special selection
of the domains (see, e.g., Miranker (1961), Dokuchaev (2007, 2010)).
\par
The paper investigates parabolic equation with homogeneous Dirichlet
boundary condition on the boundary of a domain $D\subset\R^n$ and
with mixed in time condition that connects the values  of solutions
at different times, similarly to the setting introduced in Dokuchaev
(2008) for stochastic equations.  The present paper considers a
special mixed in time conditions requiring that the solutions have a
prescribed change of profile in fixed time. Formally, this problem
does not fit the framework given by the classical theory of
well-posedness for parabolic equations based on the correct
selection of Cauchy condition. However, it is shown below that this
problem is well-posed in $L_2$-setting, and that some analog of
Maximum Principle holds. In addition, it is shown that, for any
nonnegative and non-trivial function $\g\in L_2(D)$, there exists a
unique non-negative initial function  $p(\cdot,0)$ and a number
$\a>0$ such that $p(x,0)\equiv p(x,T)+\a \g(x)$ and such that
$\int_Dp(x,0)dx=1$. This can be interpreted as an existence of a
diffusion with prescribed change of the concentration profile. An
interesting consequence is that, in the model of heat propagation, a
prescribed change of temperature during time interval $[0,T]$ can be
achieved via selection of some appropriate initial temperature, and
this problem is well-posed. On the contrary,
 a prescribed profile of temperature at time $t=T$ cannot be achieved via
selection of the initial temperature; this problem is ill-posed.
 \section{Definitions}
Let  $D\subset \R^n$ be an open bounded domain with $C^2$ - smooth
boundary $\p D$. The case when  $D$ is not connected or not simply
connected is not excluded.

Let $T>0$ be a fixed number. We consider the boundary value problems
 \baa
&\frac{\p u}{\p t}=A u, \quad &\hbox{for}\quad (x,t)\in D\times (0,T)\nonumber\\
&\hphantom{x}u(x,t)=0, \quad &\hbox{for}\quad (x,t)\in\p D\times
(0,T) \label{y4}\eaa with some additional conditions imposed at
times $t=0$ or $t=T$. Here
$$
Au\defi\sum_{i,j=1}^n a_{ij}(x,t)\frac{\p^2 u}{\p x_i \p x_j}(x,t)
+\sum_{i=1}^n f_i(x,t)\frac{\p u}{\p x_i }(x,t) -q(x,t)u(x,t).
$$ The functions $f(x,t):  D\times(0,T)  \to \R^{n}$ and $q(x,t): D\times(0,T) \to [0,+\infty)$
are measurable and bounded, such that there exist
 bounded derivatives $\p f(x,t)/\p x_i$, $i=1,...,n$. The function $a(x,t): D\times(0,T)\to
\R^{n\times n}$ is continuous, bounded, and such that there exist
 bounded derivatives $\p a(x,t)/\p x_i$, $i=1,...,n$. In
addition, we assume that the matrix $a(x,t)$ is symmetric  and
$a(x,t)\ge \d I_n$ for all $(x,t)\in D\times (0,T)$, where $\d>0$ is
a constant, and $I_n$ is the unit matrix in $\R^{n\times n}$.
\par
Problem (\ref{y4}) describes diffusion processes in domain $D$ that
are absorbed (killed) on the boundary and, with some rate, inside
$D$. The matrix $a$ represents the diffusion coefficients, the
vector $f$ describes the drift (advection), and $q$ describes the
rate of absorption inside $D$. The  assumption that $q\ge 0$ ensures
that there is absorption (loss of energy) inside the domain rather
than generation of energy.
\subsubsection*{Spaces and classes of functions}
For a Banach space $X$,
we denote  the norm by $\|\cdot\|_{ X}$.
  \par
Let $H^0\defi L_2(D)$ and  $H^1\defi \stackrel{0}{W_2^1}(D)$ be the
standard Sobolev Hilbert spaces; $H^1$ is the closure in the
${W}_2^1(D)$-norm of the set of all smooth functions $u:D\to\R$ such
that  $u|_{\p D}\equiv 0$.
\par
 Let $H^{-1}$ be the dual space to $H^{1}$, with the
norm $\| \,\cdot\,\| _{H^{-1}}$ such that if $u \in H^{0}$ then $\|
u\|_{ H^{-1}}$ is the supremum of $(u,v)_{H^0}$ over all $v \in H^1$
such that $\| v\|_{H^1} \le 1 $. $H^{-1}$ is a Hilbert space.
\par We
denote
 the Lebesgue measure and
 the $\sigma $-algebra of Lebesgue sets in $\R^n$
by $\oo\ell _{n}$ and $ {\oo{\cal B}}_{n}$, respectively.

\par
Introduce the spaces
$$
C(s,T)\defi C\left([s,T]; H^0\right),\quad \W^{1}(s,T)\defi
L^{2}\bigl([ s,T ],\oo{\cal B}_1, \oo\ell_{1};  H^{1}\bigr),
$$
and the space
$$
\V(s,T)\defi \W^{1}(s,T)\!\cap C(s,T),
$$
with the  norm $ \| u\| _{\V(s,T)} \defi \| u\| _{{\W}^1(s,T)} +\|
u\| _{C(s,T)}. $
\par
We denote the space $\V(0,T)$ as $\V$.
\begin{definition} 
\label{defsolltion2} \rm We say that equation (\ref{y4}) is
satisfied for $u\in \V$ if,  for any $t\in[0,T]$,
\baa 
\label{intur1} u(\cdot,t)=u(\cdot,0)+\int_0^t A u(\cdot,s)ds.
\eaa  The equality here is assumed to be an equality in the space
$H^{-1}$.
\end{definition}
Note that the condition on $\p D$ is satisfied in the  sense that
$u(\cdot,t)\in H^1$ for a.e.   $t$.  Further, $A u(\cdot,s)\in
H^{-1}$ for a.e. $s$. Hence the integral in (\ref{intur1}) is
defined as an element of $H^{-1}$. Therefore, Definition
\ref{defsolltion2} requires that this integral is equal to an
element of $H^0$ in the sense of equality in $H^{-1}$.
\section{The result}
\begin{theorem}
\label{Thy} For any  $\g\in L_2(D)$, there exists a unique solution
$u\in\V$ of (\ref{y4}) such that \baa\label{ppp}
u(\cdot,0)=u(\cdot,T)+\g.\eaa Moreover, there exists $c>0$  such
that \baa \|u\|_{\V}\le c\|\g\|_{L_2(D)}\label{estp} \eaa for all
$\g\in L_2(D)$.
\end{theorem}
Note that, for $u\in \V$, the value of $u(\cdot,t)$  is uniquely
defined in $L_2(D)$ given $t$, by the definitions of the
corresponding spaces. This makes  condition (\ref{ppp}) meaningful
as an equality in $L_2(D)$. By Theorem \ref{Thy}, problem
(\ref{y4}),(\ref{ppp}) is well-posed in the sense of Hadamard.
\begin{theorem}
\label{Thy2} For any non-negative and non-trivial $\g\in L_2(D)$,
the solution $u\in\V$ in Theorem \ref{Thy} is non-negative in
$D\times (0,T)$, and there exists a number $\a=\a(\g)>0$ and a
unique nonnegative solution $p\in\V$ of (\ref{y4}) such that
\baa\label{phodensity} p(\cdot,0)=p(\cdot,T)+\a \g,\qquad
\int_Dp(x,0)dx=1.\eaa  Moreover, there exists $c>0$ such that \baa
\|p\|_{\V}\le c\a\|\g\|_{L_2(D)} \label{estpa}\eaa for all $\g\in
L_2(D)$.
\end{theorem}
\par
The statement in Theorem \ref{Thy2} regarding non-negativeness of
the solution is an analog of the Maximum Principle known for
classical Dirichlet problems for parabolic equations (see, e.g.
Chapter III in Ladyzhenskaja {\it et al} (1968)).
\par
Theorem \ref{Thy2} can be applied to the model of heat propagation
in  $D$, with the loss of energy on the boundary and inside $D$ with
the rate defined by $q$. The process $p(x,t)$ can be interpreted as
the temperature at point $x\in D$ at time $t$. Therefore, Theorem
\ref{Thy2} establishes existence of the initial temperature $p(x,0)$
that ensures the prescribed change of temperature  during time
interval $[0,T]$.
\section{Proofs} For  $s\in [0,T)$ and $\xi\in
H^0$, consider the following auxiliary boundary value problem:
 \baaa
&\frac{\p v}{\p t}=A v, \quad &\hbox{for}\quad (x,t)\in D\times (s,T)\nonumber\\
&\hphantom{x}v(x,t)=0, \quad &\hbox{for}\quad (x,t)\in\p D\times
(s,T)\nonumber\\
&\hphantom{x} v(x,s)=\xi(x)&\hbox{for}\quad x\in D.\label{4.4}\eaaa
\par
Introduce  operators $\L_s:H^0\to \V(s,T)$, such that $\L_s\xi=v$,
where $v$ is the solution in $\V(s,T)$ of this problem.
\index{(ref{4.4}).} These linear operators are continuous (see,
e.g., Theorem III.4.1 in Ladyzhenskaja {\it et al} (1968)).
Introduce an operator $\Q:H^{0}\to H^0$, such that
$\Q\xi=v(\cdot,T)$, where $v=\L_0\xi$. Clearly, this operator is
linear and continuous.
\begin{lemma}
\label{lemma1}
{\rm (i)} The operator
$\Q :H^0\to H^0$ is compact; \par
{\rm (ii)} If the equation $\Q \xi=\xi$ has
the only solution $\xi=0$ in $H^0$, then the operator
$(I-\Q)^{-1}:H^0\to H^0$
is continuous.
\end{lemma}
{\bf Proof of Lemma  \ref{lemma1}}. Let $\xi\in H^0$ and $v\defi\L_0
\xi$, i.e. $v$ is the solution of the problem (\ref{4.4}). We have
that $v=\L_sv(\cdot,s)$ in $D\times (s,T)$ for all $s\in [0,T]$.
 From the second fundamental inequality
 for parabolic equations, it
follows  that \be \label{1} \|v(\cdot,T)\|_{H^1}\le
C_1\|v(\cdot,s)\|_{H^1}, \ee where $C_1$ is a positive number that
is independent from  $\xi$ and $s$ (see, e.g., Theorem IV.9.1 in
Ladyzhenskaja {\it et al} (1968)).
 Hence \baaa
\|v(\cdot,T)\|_{H^1}\le C_1\inf_{t\in[0,T]} \|v(\cdot,t)\|_{H^1} \le
\frac{C_1}{\sqrt{T}}
\left(\int_0^T\|v(\cdot,t)\|_{H^1}^2dt\right)^{1/2} \le
\frac{C_2}{\sqrt{T}}\|v\|_{\W^1(0,T)}\\ \le
\frac{C_3}{\sqrt{T}}\|\xi\|_{H^0}, \eaaa for some $C_i>0$ that are
independent from  $\xi$. Hence the operator $\Q:H^0\to H^1$ is
continuous. The embedding of $H^1$ to $H^0$  is a compact operator
(see e.g. Yosida (1965), Ch. 10.3). Then statement (i) follows.
Statement (ii) follows from Fredholm Theorem. This completes the
proof of Lemma \ref{lemma1}. $\Box$
\\
{\bf Proof of Theorem \ref{Thy}}. For $\varphi\in L_2(Q)$, consider
the problem
 \baa
&\frac{\p u}{\p t}=A u+\varphi, \quad &\hbox{for}\quad (x,t)\in D\times (s,T)\nonumber\\
&\hphantom{x}u(x,t)=0, \quad &\hbox{for}\quad (x,t)\in\p D\times
(s,T)\nonumber\\
&\hphantom{x} u(x,0)=u(x,T)&\hbox{for}\quad x\in D.\label{2004}\eaa
 By Theorem 2.2 from Dokuchaev (2004), there exists $c>0$ such that,
 for any solution $u\in\V$,
  \baaa \|u\|_{\V}\le c \|\varphi\|_{L_2(Q)} \quad \forall \varphi\in L_2(Q).
  \eaaa
Therefore, if $\g = 0$ then the only solution of (\ref{ppp}) in $\V$
is $u =0$.
 By Lemma \ref{lemma1}, it follows that
the operator $(I-\Q)^{-1}:H^0\to H^0$ is continuous.
 It follows that, for any $\g\in H^0$, there exists $\zeta=(I-\Q)^{-1}\g\in H^0$, and this $\zeta$ is
unique. Let $u\defi\L_0\zeta$. By the definitions of $\L_0$ and
$\Q$, it follows that  $u(\cdot,T)=\Q u(\cdot,0)$. We have that
$u(\cdot,0)-u(\cdot,T)=\g$, i.e.,
$$u(\cdot,0)-\Q u(\cdot,0)=\g.$$
Thus, $u\defi\L_0\zeta=\L_0(I-\Q)^{-1}\g$ is the unique solution of
(\ref{ppp}) for any $\g\in H^0=L_2(D)$. Estimate (\ref{estp})
follows from the continuity of operators $(I-\Q)^{-1}:H^0\to H^0$
and $\L_0:H^0\to \V$. The uniqueness follows from estimate
(\ref{estp}). This completes the proof of Theorem \ref{Thy}. $\Box$
\\
{\bf Proof of Theorem \ref{Thy2}}. The following definition will be
useful.
\begin{definition}
\label{def1} A function $\g :D\to\R$ is said to be piecewise
continuous if there exists a integer $N>0$ and a set of open domains
$\{D_i\}_{i=1}^N$ such that the following holds:
\begin{itemize}
\item
$\cup_{i=1}^ND_i\subseteq D\subseteq \cup_{i=1}^N\oo D_i$, and
$D_i\cap D_j=\emptyset $ for $i\neq j$. Here $\oo D_i=D_i\cup\p
D_i$.
\item For any $i\in\{1,...,N\}$, the function
$\g|_{D_i}$ is continuous and can be extended as a continuous
function $\oo \g_i: \oo D_i\cup\p D_i\to\R$.
\item For any $x\in \cup_{i=1}^N\p D_i$, there exists
$j\in\{1,...,N\}$ such that $x\in \p D_j$ and $\oo \g_j(x)=\g(x)$.
\end{itemize}
\end{definition}
\par Clearly, the set of piecewise continuous functions is everywhere dense in
$L_2(D)$, and the set of non-negative functions is closed in
$L_2(D)$. Therefore, it suffices to consider piecewise continuous
functions $\g$ only.
\par  Let $\g(x)\ge 0$ be a piecewise continuous function,
and let $u\defi\L_0(I-\Q)^{-1}\g$  be the solution of  problem
(\ref{ppp}). Since the operator $\L_0:H^0\to \V$ is continuous, we
have  that $\|u\|_{\V}\le c\|u(\cdot,0)\|_{L_2(D)}$ for some $c>0$.
It follows that if $u(\cdot,0)=0$ then $u(\cdot,T)=0$ and $\g = 0$.
By the assumptions, $\g \neq 0$. Hence $u(\cdot,0)\neq 0$ and $u
\neq 0$.
\par
Remind that  $u=\L_0\zeta$, where $\zeta=u(\cdot,0)\in H^0$. By
Theorem III.8.1 from Ladyzhenskaja {\it et al} (1968), it follows
that, for any $\e>0$, we have that $\esssup_{(x,t)\in Q'}|u(x,t)|\le
c_0$, where $Q'=\{(x,t)\in Q:\ t>\e\}$, and where $c_0>0$ depends
only on $\e, a,f,q,D$, and $\|u(\cdot,0)\|_{L_{2}(D)}$. We use here
the part of the cited theorem that deals with solutions that are
bounded on a part of the boundary; in our case, the solution
vanishes  on $\p D\times(0,T]$. It follows that \baa
\|u(\cdot,T)\|_{L_{\infty}(D)}\le c_1,\label{new}\eaa where $c_1>0$
depends only on $a,f,q,D$, and $\|u(\cdot,0)\|_{L_{2}(D)}$.
\par
Consider a sequence of functions  $u_i\in\V$ being solutions of
(\ref{y4}) such that $u_i(\cdot,0)\in C^2(\oo D)$, where $\oo D=
D\cup\p D$, $u_i|_{\p D}=0$, and that
$\|u(\cdot,0)-u_i(\cdot,0)\|_{L_2(D)}\to 0$ as $i\to +\infty$. By
Theorem IV.9.1 from Ladyzhenskaja {\it et al} (1968),
$u_i(\cdot,T)\in C(\oo D)$. (More precisely, there exists a
representative $\oo u_i(\cdot, T)$ of the corresponding element of
$H^0=L_2(D)$ which is a class of $\oo\ell_n$-equivalent functions).
By  (\ref{new}) and by the linearity od the problem, we have that
$\|u(\cdot,T)-u_i(\cdot,T)\|_{L_{\infty}(D)}\to 0$ as $i\to
+\infty$. Since the set $C(\oo D)$ is closed in $L_{\infty}(D)$,
 it follows that
there exists a  representative $\oo u$ of the corresponding element
of $\V$ such that $u(\cdot,T)$  is continuous  in $\oo D$. We have
that $u(\cdot,0)=u(\cdot,T)+\g$, hence there exists a piecewise
continuous representative of $u(\cdot,0)\in H^0$.
\par
Let us show that $\essinf_{x\in D}u(x,0)\ge 0$. Suppose that \be
\label{-} \essinf _{x\in D}u(x,0)<0. \ee If (\ref{-}) holds, then
there exists a piecewise continuous representative $\oo u(\cdot,0)$
of $u(\cdot,0)$, such  that there exists $\w x\in D$ such that
$$
\oo u(\w x,0)<0, \quad
\oo u(\w x,0)\le \oo u( x,0) \quad \hbox{for a.e.}
\ x\in D.
$$
Let $\w v\defi\oo u(\w x,0)$ be considered  as an element of
$L_2(D)$. We have that \baa \oo u(\cdot,T)=\Q \oo u(\cdot,0)=\Q \w v
+\Q (\oo u(\cdot,0)-\w v). \label{nn0}\eaa
\par
By the assumptions, $\oo u(\cdot,0)-\w v\ge 0$.  Let us show that
\baa (\Q [\oo u(\cdot,0)-\w v])(\w x)>0. \label{nn1}\eaa For this,
it suffices to show that $\oo u(\cdot,0)\neq \w v$, since a
nonnegative solution of parabolic equation (\ref{y4})  is either
identically zero or strictly positive everywhere in $D\times (0,T]$.
 Suppose that
$\oo u(x,0)\equiv \w v$. By the Maximum Principle for parabolic
equations (see, e.g. Theorem III.7.2 from Ladyzhenskaja {\it et al}
(1968)), $\oo u(x, T)=(\Q \w v)(x)\ge \w v$ for all $x$;  we apply
the version of Maximum Principle for non-positive solutions given
that $\w v<0$. It follows that if $\oo u(x,0)\equiv \w v$ then $\oo
u(x,t)$ does not satisfy (\ref{ppp}) with non-negative $\g \neq 0$.
Thus, $\oo u(\cdot,0)\neq \w v$. Hence (\ref{nn1}) holds.
\par
Further, by the Maximum Principle for non-positive solutions again,
it follows that \baa (\Q \w v)(\w x)\ge \w v=\oo u(\w x,0).
\label{nn2}\eaa By (\ref{nn0})-(\ref{nn2}), we have that $\oo u(\w
x,T)> \oo u(\w x,0). $  It follows that if (\ref{-}) holds then $\oo
u$ does not satisfy (\ref{ppp}) with $\g(x)\ge 0$. Thus, $u(x,0)\ge
0$ a.e.
\par
Let   $$\a\defi\left(\int_D u(x,0)dx\right)^{-1},\quad p\defi \a
u.$$ We have that (\ref{phodensity}) holds.
 By the linearity of problem (\ref{ppp}), it follows that
and $p(\cdot,0)-p(\cdot,T)=\a \g$ and that (\ref{y4}) holds for $p$.
Therefore, $p$ is such as required. Estimate (\ref{estpa}) follows
immediately from estimate (\ref{estp}) and from the selection $p=\a
u$.

Finally, let us show that $p$ is unique. Let $(p_i,\a_i)$ be such
that (\ref{y4}), (\ref{phodensity}) hold for $p_i\in \V$, $\a_i>0$,
$i=1,2$. Let $u_i=p_i/\a_i$. Clearly, $u_i$ is solution of
(\ref{ppp}). By the uniqueness established in Theorem \ref{Thy}, we
have that $u_1=u_2$. Hence $p_1=p_2\a_2/\a_1$. If $\a_1\neq \a_2$,
then it is not possible to have that
$\int_Dp_1(x,0)dx=\int_Dp_2(x,0)dx=1$. Therefore, $\a_1=\a_2$ and
$p_1=p_2$. This completes the proof of Theorem \ref{Thy2}. $\Box$
\subsubsection*{Acknowledgment} {This work was
supported by  The Australian Technology Network – DAAD Germany joint
research cooperation scheme.}
\section*{References} $\phantom{xxi}$Beck, J.V. (1985). {\it
Inverse Heat Conduction.} John Wiley and Sons, Inc..

\par
 Dokuchaev, N.G. (2004). Estimates for distances between first exit
times via parabolic equations in unbounded cylinders. {\it
Probability Theory and Related Fields}, {\bf 129} (2), 290 -
314.

\par Dokuchaev, N. (2007). Parabolic equations with the second
order Cauchy conditions on the boundary. {\it  Journal of Physics A:
Mathematical and Theoretical}.
 {\bf 40}, pp. 12409--12413.
\par
 Dokuchaev N. (2008). Parabolic Ito equations with mixed in time
conditions. {\it Stochastic Analysis and Applications} {\bf 26},
Iss. 3, 562--576.
\par
 Dokuchaev, N. (2010). Regularity for some backward
heat equations,   {\it Journal of Physics A: Mathematical and
Theoretical},  {\bf 43} 085201.

\par Glasko V. (1984). {\em Inverse Problems of Mathematical
physics}. American Institute of Physics. New York.

\par
Ladyzhenskaja, O.A.,  Solonnikov, V.A., and   Ural'ceva, N.N.
(1968). {\it Linear and Quasi--Linear Equations of Parabolic Type.}
Providence, R.I.: American Mathematical Society.
\par
Miller, K. (1973). Stabilized quasireversibility and other nearly
best possible methods for non-well-posed problems. {\em In:
Symposium on Non-Well-Posed Problems and Logarithmic Convexity.
Lecture Notes in Math.} V. 316, Springer-Verlag, Berlin, pp.
161--176.
\par Miranker, W.L. (1961). A well posed problem for the
backward heat equation. {\em Proc. Amer. Math. Soc.} 12 (2), pp.
243–274
\par Prilepko A.I., Orlovsky D.G.,
Vasin I.A. (1984). {\em Methods for Solving Inverse Problems in
Mathematical Physics}. Dekker, New York.
\par
Seidman, T.I. (1996). Optimal filtering for the backward heat
equation, {\em  SIAM J. Numer. Anal.} {\bf 33}, 162-170.
\par Tikhonov, A.
N. and Arsenin, V. Y. (1977). {\it Solutions of Ill-posed Problems.}
W. H. Winston, Washington, D. C.
\par
Yosida, K. (1965). {\it Functional Analysis.}  Springer, Berlin
Heidelberg New York.
\end{document}